\documentclass[11pt]{amsart} 

\usepackage{amsmath,amssymb,amsthm,amsfonts,verbatim}
\usepackage{enumerate}
\usepackage{color}
\usepackage{graphicx,url}

\theoremstyle{plain}
\newtheorem{theorem}{Theorem}[section]

\newtheorem{proposition}[theorem]{Proposition}

\newtheorem{lemma}[theorem]{Lemma}

\newtheorem{remark}[theorem]{Remark}

\newtheorem{corollary}[theorem]{Corollary}

\theoremstyle{definition}
\newtheorem{definition}[theorem]{Definition}

\newtheorem*{convention*}{Convention}

\newcommand{\dmo}{\DeclareMathOperator}

\dmo{\Mod}{Mod}
\dmo{\Ig}{\mathcal{I}_g}
\dmo{\Span}{span}
\dmo{\Diff}{Diff}
\dmo{\Homeo}{Homeo}
\dmo{\dist}{dist}
\dmo\BDiff{BDiff}
\dmo\SO{SO}
\dmo\slide{sl}
\dmo\im{im}
\dmo\id{id}
\dmo\Fix{Fix}
\dmo\Stab{Stab}
\dmo\Mcg{Mcg}
\dmo{\Hg}{\mathcal{H}_g}
\dmo{\Tg}{\mathcal{T}_g}

\renewcommand{\epsilon}{\varepsilon}

\newcommand{\ncd}{\mathcal{NC}^\dagger}

\newcommand{\cd}{\mathcal{C}^\dagger}

\newcounter{sebcomments}

\newcounter{richardcomments}

\newcounter{jonathancomments}

\title{Towards the boundary of the fine curve graph}
\author{Jonathan Bowden}
\author{Sebastian Hensel}
\author{Richard Webb}

\begin{document} 
\maketitle

\begin{abstract}
  The fine curve graph was introduced as a geometric tool to study 
  homeomorphisms of surfaces. In this paper we study the Gromov
  boundary of this space and the local topology near points associated
  with certain foliations and laminations. We then give several
  applications including finding dynamically explicit elements with
  positive stable commutator length, and proving a Tits alternative
  for subgroups of $\textrm{Homeo}(S)$ containing a pseudo-Anosov map,
  generalizing a result of Hurtado-Xue.
\end{abstract}

\section{Introduction}
\subsection*{Geometry of Curve Complexes}
The curve complex was introduced by Harvey \cite{Harvey} to study
mapping class groups, and has subsequently become one of the most
powerful tools to do so. One of the main reasons is that the coarse
geometry of the curve graph was shown to be negatively curved in
the landmark work of Masur-Minsky \cite{MM1}. More precisely, they
show that the curve graph $\mathcal{C}(S)$, i.e.\ the $1$-skeleton of
the curve complex, of any hyperbolic surface is $\delta$-hyperbolic in
the sense of Gromov \cite{Gromov}.

The curve graph (as a coarse metric object) has a natural
bordification given by considering its {\em Gromov boundary}. Furthermore, this
abstract boundary has a natural geometric description as the space of
minimal measureable (singular) foliations by Klarreich
\cite{Klarreich}.
This identification means that one can translate
coarse geometric properties of the action on the curve graph into
honest topological or dynamical properties of mapping classes acting
on the surface.

\smallskip In previous work \cite{Dagger1} we introduced the {\em fine
  curve graph} as a variant of the curve graph in order to
study surface homeomorphism and diffeomorphism groups. The fine curve
graph $\cd(S)$ is defined as the graph whose vertices are embedded
essential simple curves which are connected by an edge if they are
disjoint. A key property of this graph is that it is $\delta$-hyperbolic
and, in subsequent work with Mann and Militon \cite{dagger2}, we began to develop
a dictionary between dynamical properties of homeomorphisms
and geometric properties of their action on the fine curve graph.

\subsection*{Classification of Surface Homeomorphisms}
Thurston famously classified mapping classes \cite{Thurston} by
showing that any mapping class is isotopic to one that is
\emph{periodic}, \emph{reducible} or \emph{pseudo-Anosov}.  Masur and
Minsky \cite{MM1} re-interpreted this classification in terms of the curve graph:
pseudo-Anosov mapping classes are exactly those that act on the curve
graph \emph{hyperbolically} (i.e. with positive asymptotic translation
length). Furthermore,  the endpoints at infinity of a quasi-axis are exactly the
stable and unstable foliation of the pseudo-Anosov (under the
identification given by Klarreich's theorem). Periodic or reducible
mapping classes act \emph{elliptically} (i.e.\ with bounded diameter
orbits).  We remark that the third type of
isometries of a hyperbolic space -- \emph{parabolic} isometries, which
have unbounded orbits, but zero asymptotic translation length -- do
not occur when considering the action of the mapping class group on curve graphs.

\smallskip For the action of actual homeomorphisms on $\cd(S)$, developing such a dictionary is work in progress.
A basic tool to study the geometry of the fine curve graph is the existence of coarsely defined 1-Lipschitz projections to the {\em surviving} curve graph of a suitably punctured surface:
$$\cd(S)\longrightarrow \mathcal{C}^{surv}(S\setminus P).$$
The existence of these maps implies that homeomorphisms of $S$ which
fix a finite set of points $P$, and define a pseudo-Anosov mapping
class on $S\setminus P$ act as hyperbolic isometries on the fine curve
graph.  By now, it is known that the converse is also true (see
\cite{dagger2} for the case of the torus, and \cite{GM} for the
general case). However, in contrast to the classical setting, there
are also homeomorphisms which act as parabolic isometries on $\cd(S)$
\cite{dagger2} and although Guih\'{e}neuf-Militon \cite{Guiheneuf_Militon_parabolic} recently gave a dynamical characterisation of such isometries in the case of the $2$-torus, the situation in general is less clear.

\subsection*{Boundary points in the Fine Curve Graph}
The central goal of this paper is to prove Klarreich-like results for
the fine curve graph. While we cannot describe every boundary point,
we can at least control a large set of such points which arise from
geometric constructions. In particular, this will include the stable
and unstable foliations of pseudo-Anosov mapping classes (on $S-P$ for
any finite set $P$ for any surface $S$ of genus $\geq 1$), as well as
ending geodesic laminations.

More precisely, we show: 
\begin{theorem}[Boundary Points]\label{thm:intro-existence-bdy-pt}
  \begin{enumerate}[i)]
  \item Suppose that $q$ is a singular flat structure on $S$ where
    every cone angle is a multiple of $\pi$. Suppose that for the
    (singular) vertical foliation $\mathcal{F}$ every pseudo-leaf is
    simply-connected and contains at most one single angle-$\pi$-cone
    point. Then
    $\mathcal{F}$ determines a unique boundary point
    $\xi_{\mathcal{F}}$ in the Gromov boundary
    $\partial_\infty \cd(S)$.

    The stabilizer of $\xi_{\mathcal{F}}$ consists exactly of those
    homeomorphisms of $S$ that preserve the foliation $\mathcal{F}$ as a
    foliation of the surface.
  \item Suppose that $P \subset S$ is a finite set of points, $\rho$ a
    complete hyperbolic metric of finite area on $S - P$, and
    $\lambda \subset S-P$ a minimal geodesic lamination, which is not
    disjoint from any essential simple closed curve on $S-P$. Then
    $\lambda$ determines a unique boundary point $\xi_\lambda$ in the
    Gromov boundary $\partial_\infty \cd(S)$.

    The stabilizer of $\xi_\lambda$ consists exactly of those
    homeomorphisms of $S$ that preserve the lamination $\lambda$ as a
    subset of the surface.
  \end{enumerate}
\end{theorem}
The condition in i) that the foliation comes from a singular flat
structure is fairly mild; see Section~\ref{subsec:bsf-structures} for
a discussion. 

Also note that since the action of $\mathrm{Homeo}(S)$ on $\partial_\infty\cd(S)$ is
minimal \cite{LongTan}, Theorem~\ref{thm:intro-existence-bdy-pt} describes a dense set
of boundary points. In particular, this means that the prong structure of a sequence of elements can change under limits of boundary points corresponding to minimal measured foliations. On the other hand the situation is more rigid when considering {\em pairs} of foliations that correspond to the end points of an invariant quasi-geodesic of a hyperbolic isometry, in that the $1$-prongs under the limit must indeed converge to $1$-prongs of the limiting object; see the arguments in Section~\ref{thm:scl} for details on this argument, and applications.

\noindent A key ingredient in proving this theorem which is of
independent interest is to give a {\em geometric} interpretation of convergence of curves (seen as vertices of the fine curve graph) to these Gromov boundary points.  Describing
this precisely is somewhat technical; compare Theorem~\ref{thm:foliation-convergence} for details.

\subsection*{Stable commutator lengths} We have previously shown that there exist unbounded quasi-morphisms on the identity component $\Homeo_0(S)$ and hence there are elements with positive stable commutator length.
However, while in theory one could extract explicit elements
from our proof (and the work of Bestvina-Fujiwara \cite{BF} we rely
on), one does not obtain robust criteria to certify that specific,
dynamically meaningful homeomorphisms have positive scl.  The
following theorem remedies this, and gives such a geometric criterion.
\begin{theorem}
  Suppose that $\phi \in \mathrm{Mcg}(S\setminus P)$ is a
  pseudo-Anosov mapping class which is not conjugate to its inverse,
  and that $F:S \to S$ is a Thurston representative of $\phi$. Then
  $F$ has positive stable commutator length. 
\end{theorem}
If $F$ is a Thurston representative of a point pushing pseudo-Anosov of a
single point, then $F$ in fact has positive stable commutator length
in $\mathrm{Homeo}_0(S)$ without any further assumptions (see Remark
\ref{rem:single_push}).

We also consider the case of surfaces with boundary, thereby answering a question posed to us by Monod and Nariman:
\begin{theorem}
Let $S$ denote a connected surface with (possibly empty) boundary such that $\chi(S) <0$. Then the identity component $Homeo_0(S,\partial S)$ admits unbounded quasi-morphisms.
\end{theorem}
\noindent This more or less resolves the problem of which surface
homeomorphism groups admit quasi-morphisms, since for $S^2$, the
annulus, and the disk, one knows that they cannot exist -- by
\cite{BIP} these groups do not admit any unbounded
conjugation-invariant norms at all. For all other closed orientable
surfaces this was previously proved in \cite{Dagger1} and extended to
non-orientable surfaces with (non-orientable) genus at least three by
\cite{KK}. As a consequence of our main theorem, the case of the Klein
bottle is dealt with in \cite{Boeke}. 

\subsection*{Stabilizers of points of the Gromov boundary}
Given that boundary points can be associated to measurable foliations one can ask about the stabilizer of such points. The following extends results of Hurtado-Xue \cite[Theorem 1.10]{HurtadoXue} where $C^1$-regularity was assumed and only the case of the torus $T^2$ was considered:
\begin{theorem}
  Suppose that $G < \mathrm{Homeo}_+(S)$ is a group which contains a Thurston representative
  of a pseudo-Anosov element of $\mathrm{Mcg}(S\setminus P)$ for some (possibly empty) finite puncture set $P$. Then either $G$ contains a free group, or has an index $2$ subgroup which preserves a
  foliation on $S$.
\end{theorem}

Similarly we can show that certain boundary points cannot be stabilised by
parabolic (or hyperbolic) elements, answering a question posed to us by Yair Minsky.
\begin{proposition}
  Suppose that $\mathcal{F}$ is a foliation as in
  Theorem~\ref{thm:intro-existence-bdy-pt}~i) and $\xi_{\mathcal{F}}$
  the corresponding boundary point. Assume in addition that
  $\mathcal{F}$ is non-orientable. Then
  \begin{enumerate}[i)]
  \item $\xi_{\mathcal{F}}$ is not fixed by any parabolic isometry of
    $\cd(S)$.
  \item If $\mathcal{F}$ is not the invariant foliation of a
    pseudo-Anosov on $S\setminus P$, then $\xi_{\mathcal{F}}$ is not
    fixed by any hyperbolic isometry of $\cd(S)$.
  \end{enumerate}  
\end{proposition}

\subsection*{Acknowledgements} We would like to thank Frédéric Le Roux
for his insightful comments that greatly improved the exposition. We would also
like to thank Pierre-Antoine Guihéneuf, Katie Mann, Yair Minsky and Sam Nariman
 for stimulating questions and helpful comments, and 
thank Federica Fanoni and José Andrés Rodríguez Migueles for comments
on various drafts of this article.  The first
and second authors are supported by the Special Priority Programme SPP
2026 Geometry at Infinity funded by the DFG.  The third author was partially
supported by the EPSRC Fellowship EP/N019644/2.

\section{(Not Just) Background}
In this section, we recall some background about fine curve graphs, and singular flat structures. In both cases, we also prove a few new results which are adapted to our purposes.

\subsection{Hyperbolicity of the fine curve graph}
Let $S$ be a closed, connected, oriented surface of genus at least 1. The \emph{fine curve graph} $\mathcal{C}^\dagger(S)$ is the graph whose vertices are essential simple closed curves (up to reparametrisation), and edges correspond to disjointness (except if $S$ is a torus, in which case edges correspond to a single intersection point).

We have
\begin{theorem}[{\cite{Dagger1}}]
	The fine curve graph is Gromov hyperbolic. 
\end{theorem}
We will require a specific construction of quasigeodesics, which is implicit in previous work, but has not appeared in a concrete form. 

Namely, if $\alpha, \beta$ are two essential simple closed curves intersecting transversely, we define a \emph{bicorn (defined by $\alpha, \beta$)} to be an essential simple closed curve of the form
\[ a \cup b, \quad\quad a \subset \alpha, b \subset \beta, \]
where $a,b$ are subarcs. We emphasise that both simplicity and 
essentialness are requirements -- in general, an arbitrary choice of arcs $a,b$ will not yield a bicorn, even if they share endpoints.

Usually (in the study of non-fine curve graphs), the curves defining
bicorns are required to be transverse and in minimal position. We \emph{explicitly}
do not require the latter here -- this will be ensured by adding enough punctures when needed.

For any surface of finite type $\Sigma$ we let $\mathcal{NC}(\Sigma) \subset \mathcal{C}(\Sigma)$ denote the full subgraph whose vertices are non-separating curves. We need the following result of Rasmussen:
\begin{theorem}[{\cite{Rasmussen}}]\label{thm:rasmussen}
	There is a constant $\delta>0$ so that for any finite type surface $\Sigma$, the non-separating curve graph $\mathcal{NC}(\Sigma)$ is $\delta$--hyperbolic. In addition, if $\alpha, \beta \subset S$ are two curves in minimal position, then the set of all non-separating bicorns formed by $\alpha, \beta$ is Hausdorff close to a geodesic in the curve graph of $S$ joining $\alpha$ to $\beta$ (with uniform constants).
\end{theorem}
The last part of this theorem is not stated as such in
\cite{Rasmussen}, so we briefly explain why it follows from his
proof. Namely, Rasmussen obtains hyperbolicity of the non-separating
curve graphs using the following useful and well-known criterion (compare \cite[Theorem~3.15]{MS}, and \cite[Proposition~3.1]{Bowditch} for the precise formulation we use)
\begin{proposition}
  Let $X$ be a graph, $D>0$ a number, and suppose for each pair of
  vertices $x,y \in V(X)$ we have chosen a connected subgraph
  $\mathcal{L}(x,y)$ containing $x,y$.

  Suppose that if $d(x,y) \leq 1$, the diameter of $\mathcal{L}(x,y)$
  is at most $D$, and in addition for all $x,y,z$ we have
  \[ \mathcal{L}(x,y) \subset N_D(\mathcal{L}(x,z) \cup
    \mathcal{L}(z,y)). \] Then $X$ is $\delta(D)$--hyperbolic, and for
  each $x,y$ the subgraph $\mathcal{L}(x,y)$ is $B(D)$--Hausdorff
  close to a geodesic joining $x$ to $y$.
\end{proposition}
In \cite{Rasmussen}, Rasmussen shows that putting
$\mathcal{L}(\alpha, \beta)$ to be the graph spanned by bicorns
defined by non-separating curves $\alpha, \beta$ in minimal position,
the prerequisites of the above proposition are satisfied.  Hence,
Theorem~\ref{thm:rasmussen} follows. 

\subsection{Boundaries of Gromov hyperbolic spaces}
Here, we briefly recall the definition of the boundary at infinity of
a hyperbolic space in terms of Gromov products, which is suitable for
spaces that are not locally compact. We only very briefly recall the
required notions, and refer the reader to  \cite[III.H.3]{BH}  as well as \cite{Das_Simmons_Urbanski} for the non-proper setting.
  
\smallskip For a metric space $X$ and a basepoint $x_0 \in X$ we
recall the \emph{Gromov product}
\[ (x\cdot x')_{x_0} = \frac{1}{2}\left( d(x_0, x) + d(x_0, x') -
    d(x,x') \right). \] If $X$ is $\delta$-hyperbolic, then the Gromov
product is (up to uniform additive constants) the distance from $x_0$
to a geodesic joining $x$ and $x'$. In other words, the Gromov
product measures when geodesics from $x_0$ to $x$ and $x'$ start to diverge.

\smallskip A sequence $(a_i)$ in $X$ is called \emph{admissible} if
\[ (a_i \cdot a_j)_{x_0} \to \infty, \quad i,j \to \infty, \]
and two admissible sequences $(a_i), (b_j)$ are equivalent if
\[ (a_i \cdot b_j)_{x_0} \to \infty, \quad i,j \to \infty. \] We then
define the \emph{Gromov boundary} $\partial_\infty X$ to be the set of
admissible sequences up to equivalence. One then extends the Gromov product to boundary points by setting
$$(\xi\cdot\eta)_{x_0} = \inf_{(a_i)=\xi,(b_j)=\eta}\left(\liminf_{i,j \to \infty} (a_i\cdot b_j)_{x_0}\right)$$
One then defines a topology on $X \cup \partial_\infty X$. Namely, for
$\xi \in \partial_\infty X$, a neighbourhood basis is given by sets of
the form 
\[ U(\xi, K) = \{ \eta \in X\cup \partial_\infty X \quad
  (\eta\cdot\xi)_{x_0} > K \}. \] Note that if $X$ is proper, the
boundary defined above agrees with the more common one in terms of
geodesic rays. For non-proper spaces, it is not clear that every
Gromov boundary point (in the sense above) is the endpoint of a
geodesic ray. One could define the boundary in terms of quasi-geodesic rays, compare e.g. \cite[Proposition~2.3]{LongTan}.

\smallskip If $Y \subset X$ is any subset, we denote by $\partial_\infty Y$ the
subset of the boundary defined by accessible sequences with terms in
$Y$. If $Y$ is quasiconvex, and therefore itself hyperbolic, this is
just (an embedding of) the Gromov boundary of $Y$ (into the Gromov
boundary of $X$).

\subsection{Singular Flat and BSF Structures}\label{subsec:bsf-structures}

In this paper, we adopt the following definition, which is not the
most general possible, but fairly common and adapted for our purposes.
\begin{definition}
A \emph{singular flat structure} on a
surface $S$ will always mean a singular flat metric with discrete cone
points, all of which have cone angles that are multiples of $\pi$.  
\end{definition}
In particular, on compact surfaces, singular flat structures have
finitely many cone points.  We will often denote the set of
angle-$\pi$-singularities by $P$.

Any such structure is obtained by gluing parallel sides of a finite
number of polygons. Equivalently, such a singular flat structure can
be described by an atlas into the complex plane whose transition
functions are of the form $z \mapsto \pm z + c$, or by a (suitable)
quadratic differential on the surface.  We refer the reader to
\cite{Strebel} for background on such structures and the equivalence.

\medskip Note that in these constructions we actually obtain
slightly more than just a metric, namely the choice of a \emph{horizontal}
and \emph{vertical} direction (strictly speaking: line fields). Since this notion is crucial for the
paper, we make the following:
\begin{definition}
  A \emph{bifoliated singular flat (BSF) structure} is a tuple
  $q = (\mathcal{F}_v, \mathcal{F}_h, d)$ of a singular flat metric
  $d$, all of whose cone angles are multiples of $\pi$, and transverse
  (singular) foliations $\mathcal{F}_v, \mathcal{F}_h$ which are
  locally geodesic and orthogonal for $d$ away from the cone points.
\end{definition}
Note that the flat metric naturally endows the foliations
$\mathcal{F}_v, \mathcal{F}_h$ with transverse invariant measures,
locally induced by orthogonal distance to the leaves.

A leaf of a singular foliation is called \emph{regular} if it contains
no singular point.
A \emph{singular leaf} is a maximal leaf segment or ray bounded by
singularites without singularities in its interior. We denote by a
{\em pseudo-leaf} a connected subset of the surface which is a finite
union of singular leaves. With this
definition every nonsingular point on $S$ is contained in exactly one
regular or singular leaf (but it might be contained in several pseudo-leaves).
Since there are only finitely many singular
leaves, almost every point is regular (and almost every point on any
path transverse to the foliation is also regular).  We recall the
following definitions.
\begin{definition}
  Suppose $\mathcal{F}$ is a foliation determined by a BSF structure
  with angle-$\pi$ cone points $P$.

  The singular foliation $\mathcal{F}$ is called \emph{minimal} if
  every regular leaf is dense. It is called \emph{ending} if every pseudo-leaf
  is simply-connected and contains at most one point of $P$.
\end{definition}
Note that ending implies minimal: the boundary of the closure of an infinite, 
regular leaf is empty or contains a non-simply-connected pseudo-leaf.

The condition on pseudo-leaves only intersecting $P$ once has the
following reason: for a measured foliation, being ending is really the
condition that the isotopy class of every essential closed curve on
$S-P$ has positive transverse measure (since any such foliation
carries a transverse measure of full support). If the foliation
contains e.g. a singular leaf $\lambda$ connecting two distinct
$p_1, p_2 \in P$ then the boundary $\alpha$ of a small neighbourhood
of $s$ is an essential curve, and up to isotopy the transverse measure
can be made arbitrarily small -- but there is no representative of
transverse measure $0$ which is an embedded curves: any such sequence
degenerates to just $\lambda$.

\bigskip A natural question at this point is: which singular
foliations occur as vertical foliations of BSF structures. In the case
where there are no angle-$\pi$-singularities, this is discussed in
\cite{HubbardMasur}: every singular foliation admitting a transverse
measure of full-support is Whitehead-equivalent to one that is the
vertical foliation of a quadratic differential (and, hence, induced by
the corresponding BSF structure).  If one allows
angle-$\pi$-singularities in a set $P$, one can first pass to a
$P$-resolving cover (as in Section \ref{subsec:good}), and apply
Whitehead moves there.  We expect that all results of this paper are
also true without the requirement of being induced by a BSF structure,
but we restrict to the latter case for simplicity, and since this
class is most useful in applications.


\subsection{$P$-resolving covers and CAT(0) geometry}\label{subsec:good}
Given a BSF structure $q$ on a surface $S$, a \emph{straight segment}
is an embedded segment $I \subset S$ which intersects the
singularities of the singular flat metric at most in its endpoints,
and which is a straight line in the flat charts. A \emph{geodesic for
  $q$} is a concatenation of straight segments $s_i$ so that the angle
between $s_i$ and $s_{i+1}$ at the joining point is at least $\pi$ on
both sides. If the singular flat metric of $q$ does not have
angle-$\pi$ cone points, these are indeed locally geodesic for the
metric, and are very well-behaved (since the induced metric is locally
CAT(0)).

However, we explicity need to allow angle-$\pi$-cone points (as these
appear e.g. in the invariant foliations of point-pushing
pseudo-Anosovs) -- and then the geodesics for the singular flat metric
are much less well-behaved, since they can run into an
angle-$\pi$ singularity and return in exactly the same way.

\smallskip To control geodesics, we will therefore frequently use the
following trick.  Given a surface $S$ and a bifoliated singular flat
structure $q=(\mathcal{F}_v, \mathcal{F}_h, d)$, we consider a finite
branched cover $\hat{S} \to S$, which branches with order at least $2$
at every angle-$\pi$-cone point. We denote by
$\hat{q}=(\hat{\mathcal{F}}_v, \hat{\mathcal{F}}_h, \hat{d})$ the
lifted BSF structure. We call the surface $\hat{S}$ together with the
BSF structure $\hat{q}$ a \emph{$P$-resolving cover} of $q$. Such
$P$--resolving covers always exist (e.g. one could first take a
regular cover to make the number of angle-$\pi$ singularities even,
and then take a branched cover with branching degree $2$ at all of
them).  We emphasize that that the $P$-resolving cover is not unique,
and in the sequel we will always choose a $P$-resolving cover for any
BSF structure, although the specific choice will not matter. 

Note that, away from the cone points, the covering map is a local
isometry, and so geodesics for $q$ will lift to geodesics for
$\hat{q}$. In particular, the covering map is locally injective on
each pseudo-leaf, except at points in $P$ (here, it is important, that
the branching happens \emph{only} at $P$).

If $\mathcal{F}_v$ is ending, then $\hat{\mathcal{F}}_v$ is again
ending.  Namely, if $\hat{\mathcal{F}}_v$ had a pseudo-leaf $\lambda$
which is not simply connected or contains several angle-$\pi$
singularities, then the image of $\lambda$ would contain a problematic
pseudo-leaf for $\mathcal{F}_v$. This argument is the
reason we had to exclude paths formed by singular leaves joining
angle-$\pi$-singularities -- if such paths would exist, they could
lift to closed leaves or leaf cycles! 

The main reason we consider $P$-resolving covers is the following (compare
\cite[Theorem 14.2.1]{Strebel}).
\begin{lemma}\label{lem:cat0-geodesics}
  Suppose $q$ is a BSF structure on $S$ and $(\hat{S}, \hat{q})$ is
  a $P$-resolving cover. Then no two geodesic segments (for the flat metric
  defined by $\hat{q}$) bound a bigon.
\end{lemma}
\begin{proof}
  This follows from the fact that a singular flat metric whose cone
  angles are all of the form $\pi (k+2), k \in \mathbb{N}$ is locally
  CAT($0$), or alternately from the Gau\ss-Bonnet formula for conical
  metrics (cf.\ \cite[Theorem 14.1]{Strebel}).
\end{proof}
This has the following useful corollary
\begin{corollary}\label{cor:common-segment}
  Suppose that $\gamma_1, \gamma_2$ are closed geodesics on
  a $P$-resolving cover $\hat{S}$. Let $\widetilde{S}$ be the universal cover of $\hat{S}$,
  and $\widetilde{\gamma}_i$ be lifts of the $\gamma_i$. Then the intersection $ \widetilde{\gamma}_1 \cap \widetilde{\gamma}_2 $ consists of a common geodesic segment.

  Furthermore, between any two points in $\widetilde{S}$ there is a
  unique geodesic for the singular flat metric (which realises the distance for the path metric).
\end{corollary}
\begin{proof}
  If the intersection were disconnected, we could find a bigon bounded
  by subsegments of the $\widetilde{\gamma}_i$. Similarly, if
  geodesics were not unique, we could find a bigon bounded by two
  distinct geodesics joining the same endpoints. Existence of
  geodesics follows from completeness of the metric by a standard
  Arzela-Ascoli argument.
\end{proof}
Note that on (the universal cover of) $S$ itself, these claims are false
if $q$ has angle--$\pi$ cone points -- geodesics can bound bigons
containing angle-$\pi$-cone points -- hence the need for $P$-resolving
covers.

\begin{corollary}\label{cor:bi-infinite}
  Let $\widetilde{S}$ be the universal cover of the $P$-resolving cover
  $\hat{S}$. Then the metric on $\widetilde{S}$ induced
  by the lifted BSF structure is Gromov hyperbolic.
  
  Suppose that $\gamma$ is a bi-infinite geodesic for this metric.  
  Then $\gamma$ separates $\widetilde{S}$ into two disks, each of
  which is convex.
\end{corollary}
\begin{proof}
  Since $\hat{S}$ is a compact surface, the metric on $\widetilde{S}$ induced
  by the lifted singular flat structure is quasi-isometric to the hyperbolic
  metric on the disk. 

  The separation property is inherited from the corresponding fact for
  hyperbolic quasigeodesics, since geodesics in the flat metric are quasi-geodesics for the hyperbolic metric. The convexity claim then follows from
  Corollary~\ref{cor:common-segment}.
\end{proof}

These results in particular imply that leaf spaces of the horizontal
and vertical foliations of $\widetilde{q}$ equipped with the
transverse measures are real trees $T_h, T_v$. We denote by
\[ \pi_h : \widetilde{S} \to T_h, \quad \pi_v : \widetilde{S} \to
  T_v \] the corresponding projections. Since the trees are
\emph{dual} to the laminations, the diameter of the (vertical)
projection $\pi_v(X)$ of the projection of a set $X$ measures the
``horizontal width'' (i.e. how much the set crosses the vertical
leaves).  We will not need any specific properties of these trees
and refer the interested reader to \cite{MorganShalen} for
further details. For us, the projections to the trees are just a convenient
way to measure the size of sets lifted to the universal cover relative to the foliations.

\subsection{Hyperbolic elements and quasi-invertibility}\label{sec:BF}
We recall the notion of coarse equivalence introduced by Bestvina-Fujiwara \cite{BF}. For this consider hyperbolic isometries $f,g$ of a general $\delta$-hyperbolic space $(X,d)$. These then have invariant quasi-axes given for example by taking the orbit of an arbitrary point say $A= (\alpha_i), A'=(\beta_j)$ of quality $L,L'$ respectively. Then we say $f \sim g$ are {\em quasi-equivalent} if  there is a (fixed) constant $B= B(L,L',\delta)$ such that for any $D > 0$ there is an isometry moving a segment of length $D$ in $A$ within a $B$-neighbourhood of $A'$, respecting the direction of $A$ and $A'$. One can check, using the hyperbolicity condition, that this definition is well-defined and indeed gives an equivalence relation.

Since quasi-geodesics determine points on the Gromov boundary, we can set $\xi_{\pm} = \lim_{i \to \pm \infty} \alpha_i$ and  $\eta_{\pm} = \lim_{j \to \pm \infty} \beta_j$. This definition is then equivalent to the following: There exist $g_k \in \textrm{Isom}(X,d)$ so that
$$\lim_{k \to \infty} g_k (\xi_{\pm}) = \eta_{\pm}.$$
In the special case $f\sim f^{-1}$ we say that $f$ is {\em quasi-invertible}.

The following is a consequence of the fact that the action of the mapping class group on the curve graph satisfies the WPD condition and can be deduced from \cite[Propositions 6, 11]{BF}:
\begin{theorem}\label{thm:BF1}
Let $G \subset Mcg(\Sigma)$ be a subgroup of the mapping class group of a finite type surface $\Sigma$ that is not virtually cyclic and contains a pseudo-Anosov element. Then $G$ contains non-quasi-invertible pseudo-Anosovs. 
\end{theorem}
\noindent Below we will need the following, which states that if one can exclude quasi-invertibility then one can construct quasi-morphisms and certify positive stable commutator length.
\begin{theorem}[\cite{BF} Proposition 5]\label{thm:BF2}
Let $G \curvearrowright (X,d)$ be an action by isometries on a $\delta$-hyperbolic space. If $f \not\sim f^{-1}$ is a non-quasi-invertible isometry, then there is a quasi-morphism that is unbounded on the group generated by $f$. In particular, $f$ has positive stable commutator length.
\end{theorem}

\section{The Target Lemma}\label{sec:target}
In this section, we fix a BSF structure $q$ on a surface $S$. We choose
$(\hat{S},\hat{q})$ a $P$-resolving cover and let $\widetilde{q}$ the
corresponding lifted singular flat structure on the universal cover of
$\hat{S}$. We assume throughout that the vertical foliation of $q$
(and therefore $\hat{q}$) is ending, and so vertical geodesic rays
that do not have compact closure intersect every open set. The goal of
this section is to prove a version of quantitative minimality: if $X$
is a set which is ``horizontally wide, but vertically small'', then a
vertical geodesic (parametrized by arclength) hits $X$ after a bounded
time.

The following notions will make this precise:
\begin{definition}\label{def:size_width}
  Let $\gamma$ be a (continuous) path on $S$, and let $\widetilde{\gamma}$
  be a lift of $\gamma$ to the universal cover $\widetilde{S}$.
  \begin{enumerate}
  \item The \emph{total ($q$-)width} (respectively \emph{total
      ($q$-)height}) is the diameter of the projection
    $\pi_v(\widetilde{\gamma}) \subset T_v$ in the vertical dual tree
    (respectively the diameter $\pi_h(\widetilde{\gamma}) \subset T_h$
    in the horizontal dual tree). 
  \item The \emph{($q$-)size} of $\gamma$ is the diameter of
    $\widetilde{\gamma}$ for the metric defined by $\widetilde{q}$.
  \item If $\alpha$ is a simple closed curve, we define the
    \emph{size} and \emph{width of $\alpha$} as the
    maximal size or width when
    interpreting $\alpha$ as a path (i.e. by choosing a basepoint).
  \end{enumerate}
  We emphasize that in (3), different choices of basepoint can change
  the size and width of the loop representing $\alpha$ by at most a
  factor of $2$ -- for our argument the specific choice will never
  play a significant role, and we could just as well have chosen the minimum.
\end{definition}
\begin{remark}
  We note that the notion of size coarsely does not depend on the flat
  metric (as long as the angle-$\pi$ singularity set $P$ in
  unchanged), since it is essentially a topological notion: one could
  instead choose a fundamental domain in the universal cover of $S-P$
  and consider the diameter of the graph of those fundamental domains
  which are crossed by a lift. For us, the metric point of view is
  more convenient though.
\end{remark}
\begin{remark}
  Given a $P$-resolving cover $(\hat{S},\hat{q})$, width and size are clearly
  well-defined since the deck group of the universal cover of
  $\hat{S}$ acts by isometries, preserving the horizontal and vertical
  foliations and their transverse measures.

  Choosing a different $P$-resolving cover can only change these
  quantities by a multiplicative constant depending on the degree
  of the $P$-resolving cover (which is bounded). For us, the specific values are rarely relevant,
  so the explicit choice of $P$-resolving cover will not have any influence.
\end{remark}
\begin{remark}
  Note that size and width are well-behaved with respect to smoothing,
  and tubular neighbourhoods. Namely, if $\gamma'$ is contained in an
  $\epsilon$--neighbourhood of $\gamma$ (with respect to the flat
  metric), then both total width and size of $\gamma'$ are at most
  $2\epsilon$ larger than the corresponding quantity for $\gamma$.
\end{remark}
The following is a first instance of the quantitative minimality
mentioned in the beginning of this section.
\begin{lemma}[Fast Return Lemma]\label{lem:fast-return}
  Let $q$ be a BSF structure whose vertical foliation is
  ending. Then for any $B,\epsilon>0$ there is an $L>0$ so that the
  following holds: Suppose $g \subset S$ is a geodesic of length at
  most $B$, and width at least $\epsilon$. Then any vertical unit
  speed geodesic of length at least $L$ intersects $g$.
\end{lemma}
\begin{proof}
  This is a compactness argument: if such a constant does not exist,
  one could find a sequence of geodesics $g_n$ of length at most $B$
  and width at least $\epsilon$, and vertical geodesics $\sigma_n$ of
  length $n\to \infty$ which are disjoint from $g_n$. By Arzela-Ascoli
  we can take a (sub)limit, which converges to a nontrivial geodesic
  $g_\infty$ having width at least $\epsilon$, and a vertical
  (possibly singular) ray $\sigma_\infty$ which does not intersect
  $g_\infty$. This contradicts minimality of the vertical foliation.
\end{proof}
\begin{lemma}[Target Lemma]\label{lem:target-lemma}
  Let $q$ be a BSF structure on $S$ with ending vertical
  foliation. Then, for any $B, \epsilon$ there is a number $L>0$ with
  the following property. Suppose that $\gamma:[0, 1] \to S$ is a path
  so that the total width of $\gamma$ is at least $\epsilon$ and the
  size of $\gamma$ is at most $B$.
  
  Then any vertical geodesic for $q$ of length at least $L$ intersects
  $\mathrm{im}(\gamma)$.
\end{lemma}

\begin{remark}
  Observe that the Target Lemma~\ref{lem:target-lemma} is (in some
  sense) optimal: by starting with any transverse arc of total width
  $\epsilon$ and pushing part of it to follow the vertical foliation
  (pushing more subarcs out of the way if intersections are created)
  we can generate an arc where \textbf{some} first return time is
  enormous, without changing the total width (but making the size
  enormous at the same time).
\end{remark}

It is likely that the target lemma could be proved by a contradiction
argument directly, but we instead opt to reduce it to the Fast Return
Lemma above.
\begin{proof}
  Denote by $(\hat{S}, \hat{q})$ the chosen $P$-resolving cover of
  $(S,q)$. Since we defined width and size using the $P$-resolving cover, it
  suffices to show the conclusion of the lemma for
  $(\hat{S}, \hat{q})$ itself.
  
  In other words, we may assume that the metric of $q$ has no
  angle-$\pi$ cone points, so that
  Corollaries~\ref{cor:common-segment} and \ref{cor:bi-infinite}
  apply.

  \smallskip Lift $\gamma$ to a path $\widetilde{\gamma}$ in the
  universal cover. Denote by $l_0, l_1$ the (possibly singular)
  vertical leaves through the endpoints of $\widetilde{\gamma}$ which
  are limits of regular leaves through points of $\widetilde{\gamma}$.
  Note that these are disjoint bi-infinite geodesics, each of which
  separates the disk (Corollary~\ref{cor:bi-infinite}).
  
  By possibly replacing $\gamma$ with a subpath, we may assume that 
  \begin{enumerate}
  \item $\widetilde{\gamma}$ intersects $l_0, l_1$ only in its
    endpoints, and
  \item the distance between $l_0, l_1$ in $T_v$ is at least
    $\epsilon$.
  \end{enumerate}
  
  Denote by $\sigma$ the shortest geodesic joining $l_0$ and $l_1$ and
  observe that by (2) it therefore also has total width
  $\geq \epsilon$ (note that it might be the case that any geodesic
  joining $l_1, l_0$ is singular, and necessarily not just
  horizontal). On the other hand, note that the length of $\sigma$ is
  at most $B$, since $\widetilde{\gamma}$ also joins $l_0, l_1$, and
  the latter is a diameter $B$ subset by assumption.
  
  \smallskip The Gromov hyperbolicity of the singular flat metric on
  $\widetilde{S}$ implies that some point on $\widetilde{\gamma}$ has
  distance at most $K_1=K_1(B,\epsilon)$ from $\sigma$ (as any path
  \emph{outside} a $k$--neighbourhood of the shortest connection
  between the two geodesics has diameter bounded from below by an
  exponential function of $k$). Thus, the distance between any point
  on $\widetilde{\gamma}$ and any point on $\sigma$ is at most
  $K_2 = K_1 + 2B$.

  \smallskip Next note that for any point $p \in \sigma$, one of the
  vertical half-leaves through $p$ intersects $\widetilde{\gamma}$
  (since both $\widetilde{\gamma}$ and $\sigma$ join the vertical
  geodesics $l_0,l_1$). Since geodesics in $\widetilde{S}$ are unique
  and realise distance (Corollary~\ref{cor:common-segment}), there is
  in fact a vertical segment of length at most $K_2$ starting in $p$
  which intersects $\widetilde{\gamma}$. This shows that any vertical
  segment starting in a point of $p\in \sigma$ which extends length
  $K_2$ in both directions needs to intersect $\widetilde{\gamma}$.
  
  \smallskip Since $\sigma$ has length at most $B$ and width at least
  $\epsilon$, the Fast Return Lemma~\ref{lem:fast-return} applies to
  $\sigma$, and guarantees that any vertical geodesic of length
  $\geq L$ intersects $\sigma$.

  Combining the previous two observations, we conclude that any
  vertical geodesic of length at least $2L+2K_2$ intersects $\gamma$,
  proving the lemma.
\end{proof}

We also need a version of the Target Lemma for laminations.  Namely,
let $P \subset S$ be a finite set of points, and $\rho$ a complete
hyperbolic metric of finite area. We say that a measured geodesic
lamination $\lambda$ on $S-P$ is \emph{ending relative to $P$} if
every leaf is dense, and additionally no essential simple closed curve
of $S-P$ is disjoint from $\lambda$. Note that \emph{essential} means
that the curve does not bound a disk or once-punctured disk in $S\setminus P$.

The requirement to admit a transverse measure (of full support) is no
restriction: any minimal geodesic lamination carries a transverse measure.

To define a \emph{($P$-)resolving cover} in this setting, we take a branched
cover at every point in $P$ of degree $3$. The lifted geodesic
lamination $\hat{\lambda}$ has the property that every complementary
region is a once-punctured geodesic $k$--gon with $k \geq 3$.
We also choose a hyperbolic metric $\rho_0$ of finite area on
$\hat{S}$ \emph{without cusps}. We then use the universal cover of
$\hat{S}$ with the lifted hyperbolic metric $\widetilde{\rho}_0$ to
define \emph{size}, and the transverse measure of the lift
$\widetilde{\lambda}$ of $\hat{\lambda}$ to define \emph{width}. Note
that the lamination $\hat{\lambda}$ is of course not geodesic for
$\rho_0$, but it is still quasigeodesic:
\begin{lemma}
  Let $\widetilde{S}$ be the universal cover of $\hat{S}$ and
  $\widetilde{\rho}_0$ the lift of the hyperbolic metric
  $\rho_0$. Then every leaf of the lift of $\hat{\lambda}$ to
  $\widetilde{S}$ is a quasi-geodesic for $\widetilde{\rho}_0$.
\end{lemma}
\begin{proof}
  One way to see this is to complete $\hat{\lambda}$ to a singular
  foliation $\mathcal{F}$, where every singularity has at least $3$
  prongs.  Theorem~2.3 of \cite{HubbardMasur} shows that $\mathcal{F}$
  is Whitehead equivalent to a foliation $\mathcal{F}'$ which is the
  vertical foliation of a BSF structure (in their language: the
  vertical foliation of a holomorphic quadratic differential). Arguing
  as before, the leaves of $\mathcal{F}'$ therefore lift to
  quasi-geodesics in $\widetilde{S}$. The Whitehead moves are
  supported on disks in $\hat{S}$, and therefore undoing them does not
  change the property of being quasigeodesic.
\end{proof}
The proof of Lemma~\ref{lem:target-lemma} now works \emph{mutatis
  mutandis}, and shows:
\begin{lemma}[Target Lemma for Laminations]\label{lem:target-lamination}
  Suppose that $\lambda$ is ending relative to $P$. Then, for any
  $B, \epsilon$ there is a number $L>0$ with the following
  property. Suppose that $\gamma:[0,1] \to S$ is a path so that the
  total width of $\gamma$ is at least $\epsilon$ and the size of
  $\gamma$ is at most $B$. Then any leaf segment of $\lambda$ of
  length at least $L$ intersects $\gamma$.
\end{lemma}

\section{Constructing Boundary Points from Foliations and Laminations}
\label{sec:construct}
In this section we work with the \emph{non-separating fine curve graph}
$\mathcal{NC}^\dagger(S)$, which is the full subgraph whose vertex set
consists of non-separating simple closed curves. This is purely for
convenience, and to be able to cite certain results from the
literature. The main results will remain true for the full fine curve
graph as well, as the two are quasi-isometric.

We will first show that a vertical foliation $\mathcal{F} = \mathcal{F}_v$ of a Bifoliated Singular Flat (BSF) structure that is {\em ending} naturally defines a unique Gromov boundary
point of the fine curve graph. The case of ending laminations will then follow {\em mutatis mutandis}.
The main result for foliations is as follows:
\begin{theorem}\label{thm:foliation-convergence}
  Suppose $q=(\mathcal{F}_v, \mathcal{F}_h, d)$ is a BSF structure
  with ending $\mathcal{F}_v$. Then there is a point
  $\xi_{\mathcal{F}_v}$ on the Gromov boundary of $\mathcal{NC}^\dagger(S)$ with the following properties:
  \begin{enumerate}
  \item $\xi_{\mathcal{F}_v}$ depends only on the foliation $\mathcal{F}_v$,
  \item a homeomorphism fixes $\xi_{\mathcal{F}_v}$ as a point on the Gromov boundary if and only if it preserves the foliation $\mathcal{F}_v$,
  \item a sequence of curves $(\beta_i)$ converges to
    $\xi_{\mathcal{F}_v}$ (as vertices of the non-separating fine curve graph) if and
    only if the size of $\beta_i$ diverges to infinity, and for any
    $B, \epsilon$ there is a $N$ so that if $b \subset \beta_i, i > N$
    is a segment of size at most $B$, then it $\epsilon$--fellow
    travels a leaf segment of $\mathcal{F}_v$ with respect to the
    metric $d$.
  \end{enumerate}
\end{theorem}

In order to construct this boundary point, we will use (quasi)-minimising properties of bicorn paths:
\begin{lemma}\label{lem:bicorns}
  Let $P \subset S$ be a finite set. Suppose that $\alpha, \alpha'$
  are two non-separating curves in minimal position on the punctured
  surface $S\setminus P$. Then any sequence of bicorn surgeries formed
  by $\alpha, \alpha'$ define a uniform (unparametrised) quasigeodesic
  in $\mathcal{NC}^\dagger(S)$. This means: there is a constant $K>0$
  so that if $(a_i)$ is a sequence of bicorn surgeries of $\alpha$
  towards $\alpha'$, then the sequence admits a reparametrization
  $i(n)$, so that $n \to a_{i(n)}$ is a
  $K$--quasigeodesic.
\end{lemma}
We also note that any bicorn formed by
$\alpha, \alpha'$ does appear in a sequence as in Lemma~\ref{lem:bicorns}.
\begin{proof}
  We may add points to $P$ so that each complementary component of
  $\alpha \cup \alpha'$ contains at least one point of $P$.
  
  Denote by $\gamma_i$ the bicorn path of non-separating curves, and
  $i(n)$ the parametrization so that they define a $K$--quasigeodesic
  for the constant from Lemma~\ref{lem:bicorns}. Given two
  indices $n,m$, the bicorns $\gamma_{i(n)}, \gamma_{i(m)}$ are not in
  minimal position on $S\setminus P$ since they might have segments in
  common. However, slight pushoffs $\gamma_n', \gamma_m'$ (which are
  distance at most $2$ in $\mathcal{NC}^\dagger(S)$) are in minimal
  position by our assumption on $P$, and still define the same isotopy
  classes in $S\setminus P$. Therefore, we have by \cite{Dagger1} that
  the distance in $\mathcal{NC}^\dagger(S)$ of $\gamma_n', \gamma_m'$
  (and hence, up to an additive error of $4$, the distance of
  $\gamma_{i(n)}, \gamma_{i(m)}$) is equal to the distance of
  $[\gamma_{i(n)}], [\gamma_{i(m)}]$ in $\mathcal{NC}(S\setminus P)$.
  The latter is then quasi-comparable to $|n-m|$ by Theorem~\ref{thm:rasmussen}.
\end{proof}

Given any set $L \subset S$, define $\mathcal{D}(L)$
to be the set of curves disjoint from $L$:
\[ \mathcal{D}(L) = \{ \alpha \in \mathcal{NC}^\dagger(S): \alpha \cap L
  = \emptyset \}. \]

\begin{lemma}\label{lem:disjoint-quasiconvex}
  The set $\mathcal{D}(L)$ is $K$--quasiconvex for some constant $K$
  independent of $L$.
\end{lemma}
\begin{proof}
  The set $\mathcal{D}(L)$ is clearly invariant under bicorn
  surgery. Hence, the result follows from Lemma~\ref{lem:bicorns}.
\end{proof}
The following then shows point (1) in Theorem \ref{thm:foliation-convergence}.
\begin{proposition}\label{prop:boundary-from-foliation}
Consider an ending foliation $\mathcal{F}= \mathcal{F}_v$ coming from a BSF-structure. Then there is a unique point of the Gromov boundary
  $\xi_\mathcal{F} \in \partial_\infty \mathcal{NC}^\dagger(S)$ contained in
  every $\partial_\infty \mathcal{D}(L)$ where $L$ is a finite union of leaf
  segments of $F$.

  In fact, if $L_i$ is a sequence of  leaf
  segments of unboundedly increasing length, then any sequence $\alpha_i \in \mathcal{D}(L_i)$ is an
  admissible sequence defining the boundary point $\xi_\mathcal{F}$.
\end{proposition}
\begin{proof}
  Let $\alpha_0$ be a basepoint in $\mathcal{NC}^\dagger(S)$. We may
  assume that $\alpha_0$ is transverse to
  $\mathcal{F}$. Consider an increasing sequence
  $L_1 \subset L_2 \subset \cdots$ of leaf segments with diverging
  length.

  We first claim that the distance of $\mathcal{D}(L_i)$ to the
  basepoint $\alpha_0$ diverges. Namely, suppose that $X \to S$ is a
  given finite (unbranched!) cover. Since $\mathcal{F}$ is ending,
  hence minimal, and $X$ is a finite cover, there is an index $n$ so
  that there is a lift $\overline{L}_n$ of $L_n$ to $X$ which is
  filling with any lift $\overline{\alpha}_0$ of $\alpha_0$ -- that
  is: every complementary component of
  $\overline{\alpha}_0\cup \overline{L}_n$ is simply-connected.

  If $\overline{\delta}$ is the lift of a curve
  $\delta \in \mathcal{D}(L_n)$ to $X$, then it is an essential curve
  which is disjoint from $\overline{L}_n$. Hence, it cannot be
  disjoint from any lift $\overline{\alpha}_0$ of $\alpha_0$.  By the
  covering criterion \cite[Lemma 4.3]{dagger2} the claim follows. 

  \smallskip Now, take $\alpha_i$ to be a non-separating,  simple curve disjoint from $L_i$.
  Since $\mathcal{D}(L_i)$ is quasiconvex by
  Lemma~\ref{lem:disjoint-quasiconvex}, the distance of
  $\mathcal{D}(L)$ to a basepoint gives a coarse lower bound for the
  Gromov product of any two points in $\mathcal{D}(L)$. Hence, 
  $(\alpha_i)$ gives an admissible sequence converging to a boundary
  point $\xi$. By the same argument, a different choice of $\alpha'_i$
  gives an equivalent sequence, hence the same boundary point. This
  shows that $\xi_\mathcal{F}$ is the only Gromov boundary point contained in the
  boundary of all $\mathcal{D}(L_i)$.
  
  \smallskip To finish the proof, it suffices to observe that if $L$
  is any (possibly disconnected) finite union of leaf segments, we can take a sequence consisting of non-separating curves $\alpha_i$ that are disjoint from $L \cup L_i$.
\end{proof}

The previous proposition is not optimal, and we now give a better
criterion for convergence.
\begin{corollary}\label{cor:convergence-crit}
  Suppose $\xi_\mathcal{F}$ is the boundary point defined by the vertical
  foliation $\mathcal{F}$ of a BSF-structure that is ending. Let $(\alpha_i)$ be a sequence of curves with size diverging to
  infinity, and assume that for all $B, \epsilon > 0$ there is an 
   $N=N(B,\epsilon)$ so that any subsegment of
  $\alpha_i, i > N$ of size $\leq B$ has width at most $\epsilon$. 

  Then $\alpha_i$ converges to $\xi_\mathcal{F}$ in the Gromov boundary.
\end{corollary}
\begin{proof}
  Let $\mathcal{D}(\epsilon, B)$ be the set of all non-separating closed curves on $S$ with the
  property that $\beta \in \mathcal{D}(\epsilon, B)$ if any
  size-$\leq B$-subsegment of $\beta$ has width at most
  $\epsilon$. 

  Since there are only finitely many isotopy classes of simple closed
  curves of size at most $B$, and $\mathcal{F}$ is ending, there is an
  $\epsilon_B$ so that any nonseparating curve of size at most $B$ has
  width at least $\epsilon_B$. For otherwise, up to taking subsequences, we would find a sequence of essential closed curves that converge in the Hausdorff sense into a (pseudo-)leaf, giving a closed essential (pseudo-)leaf of $\mathcal{F}$. In
  particular, if $\epsilon < \epsilon_B$, all curves in
  $\mathcal{D}(\epsilon, B)$ have size at least $B$.
  
  Note that a size-$B$-segment with width at most $\epsilon$ is
  contained in an $\epsilon$--neighbourhood of a length-$B$ pseudo-leaf
  segment of $\mathcal{F}$.  Hence, arguing as in the first part of
  the proof of Proposition~\ref{prop:boundary-from-foliation}, we see
  that the distance of the sets $\mathcal{D}(\epsilon, B)$ to a given
  basepoint in the fine curve graph diverges, as
  $\epsilon \to 0, B \to \infty$.
  
  Finally, suppose that
  $\beta_1, \beta_2 \in \mathcal{D}(\epsilon, B)$, and $\gamma$ is a
  bicorn formed by $\beta_1, \beta_2$. Any size-$\leq B$ subsegment of
  $\gamma$ is the union of at most two size-$\leq B$ subsegments of
  $\beta_1, \beta_2$. Thus,
  \[ \gamma \in \mathcal{D}(2\epsilon,B). \] In particular, if
  $2\epsilon < \epsilon_B$, the bicorn itself has size at least $B$,
  and thus if $B$ is large enough and $\epsilon$ is small enough, then the Gromov
  product of $\beta_1, \beta_2$ with respect to some base point is very large.

  On the other hand, given any vertical leaf segment $L$ and any
  $\epsilon, B$, we can find a curve
  $\alpha \in \mathcal{D}(L) \cap \mathcal{D}(\epsilon, B)$ by taking
  a very long vertical leaf segment disjoint from $L$ and closing with
  horizontal segment of size at most $\epsilon$ (which is possible by
  minimality).

  Together this shows that $\xi_\mathcal{F}$ is the unique Gromov boundary point
  contained in the boundary of all $\mathcal{D}(\epsilon, B)$, showing
  the corollary.
\end{proof}

We next aim to describe the converse: that is we wish to characterise what it means that a sequence of
curves converges to a boundary point coming from an ending foliation. From now on we also
choose, once and for all, a basepoint $\beta_0$ in $\ncd(S)$.

The key result is the following lemma. Both in the lemma, and its
consequences, we will assume that all curves are pairwise in general position. This can either be achieved by first smoothing and applying smooth transversality or applying topological transversality directly. Since these modifications can be achieved by a $C^0$-small change, which in turn corresponds to a change of distance at most $2$ in the fine curve graph, this is immaterial for all arguments below.
\begin{lemma}[Small Width Lemma]\label{lem:fellowtravel}
  Consider an ending foliation $ \mathcal{F}_v$ coming from a
  BSF-structure.  For any $\epsilon, B$ there is a $K$ so that the
  following holds. If $\beta$ is a simple closed curve and the Gromov
  product satisfies $(\beta \cdot \mathcal{F}_v)_{\beta_0} > K$, then every
  segment of $\beta$ of size $B$ has total horizontal width at most
  $\epsilon$, and thus $\epsilon$--fellow travels a (possibly
  singular) leaf of $\mathcal{F}_v$.
\end{lemma}

To prove this, we need a few ingredients. 
\begin{lemma}\label{lem:size-bounds-distance}
  Suppose that $q$ is a BSF structure on $S$. Let $\alpha_0$
  be a non-separating curve. Then for any $L> 0$ there is an $R>0$ so that the
  following is true: if $\beta$ is a curve of size at most $L$, then
  \[ d_{\cd}(\beta, \alpha_0) \leq R. \]
\end{lemma}
The idea of this proof is that for a curve $\beta$ of bounded size, we
can remove bigons with $\alpha_0$ via a uniformly bounded number of surgeries, and for a
configuration without bigons the number of intersections can be bounded in terms of its size.
\begin{figure}
  \centering
  \includegraphics[width=0.6\textwidth]{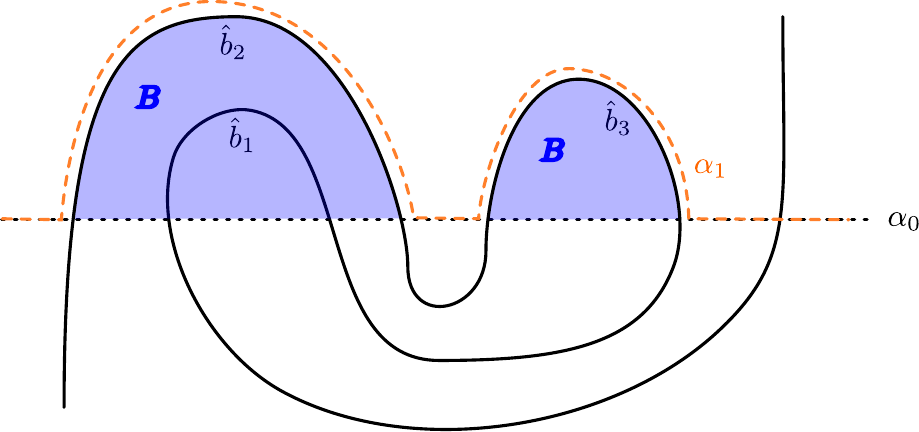}
  \caption{Removing bigons in the proof of Lemma~\ref{lem:size-bounds-distance}}
  \label{fig:bigons}
\end{figure}
\begin{proof} 
 First observe that we may assume that $\beta$ is in general position with respect to $\alpha_0$.
  Now, let $b$ be a path which lifts $\beta$ with respect to a base point not lying on $\alpha_0$, and observe that it has
  diameter at most $L$. A ball of diameter $L$ can intersect at most
  $K$ distinct lifts of $\alpha_0$ (where $K_L$ depends only on $L$ and
  $\alpha_0$).

  Let $a_1, \ldots, a_k$ (with $k \leq K_L$) be the distinct lifts
  of $\alpha_0$ which $b$ intersects. Suppose that there is some $a_i$ which is intersected more than once
  by $b$, and let $a_i$ be an outermost such lift. Let $b_1, \ldots, b_M$ denote
  all subarcs of $b$ which lie on the outer side of $a_i$ (i.e. the
  side not containing the initial point of $b$).

  Consider the images $\hat{b}_1, \ldots, \hat{b}_M$ of the $b_i$ on
  the surface $S$. These define arcs of $\beta$ which all end on the
  same side of $\alpha_0$, and do not intersect $\alpha_0$ in their
  interior.

  Thus, there is a disjoint union $B$ of (maximal) bigons based at
  $\alpha_0$ bounded by a subset of the $\hat{b}_i$ and subarcs of
  $\alpha_0$, so that all $\hat{b}_i$ are contained in $B$ (compare Figure~\ref{fig:bigons}).

  Let $\alpha_1$ be the result of surgering $\alpha_0$ at $B$, and
  pushing slightly off $B$. We observe the following points:
  \begin{enumerate}
  \item The distance between $\alpha_1$ and $\alpha_0$ in the fine
    curve graph is at most $2$,
  \item if a subarc of $\beta$ intersects $\alpha_1$, it also
    intersects $\alpha_0$ in the same points, and
  \item the subarcs $\hat{b}_1, \ldots, \hat{b}_M$ do not intersect
    $\alpha_1$.
  \end{enumerate}
  As a consequence, the number of lifts of $\alpha_1$ which $b$
  intersects is also at most $k$, and in fact the number of such lifts
  which are intersected more than once has decreased by at least one.

  Hence, after repeating this process at most $k$ times, we find a
  curve $\alpha_k$ (of distance at most $2k$ to $\alpha_0$) with the
  property that $b$ intersects at most $k\le K_L$ lifts of $\alpha_k$, and
  each in at most $1$ point.

  Hence, $\beta$ and $\alpha_k$ intersect in at most $K_L$ points, and
  therefore the distance in the fine curve graph is at most $3K_L+1=:R$.
\end{proof}

\begin{lemma}\label{lem:bicorn-existence}
  Assume that $\mathcal{F}= \mathcal{F}_v$ is the vertical foliation of a BSF-structure that is ending.
  Suppose that $\beta_0$ is a compact embedded segment which is not
  vertical. Then for some $p\in \beta_0$, the curve formed by a
  subarc of $\beta_0$ and a first return vertical flow-line starting
  in $p$ is non-separating.
\end{lemma}
\begin{proof}
  Denote by $[l_1], \ldots, [l_k]$ the isotopy classes as maps to the pair  $(S,\beta_0$) of arcs formed by
  first-return vertical flow-lines to $\beta_0$. Since they are all
  disjoint, it is clear that there is only a finite number of
  them. Collapsing $\beta_0$ to a point the $l_1, \ldots, l_k$ give an
  embedded graph on the surface (with a single vertex). By the
  ending assumption on the vertical foliation all complementary
  regions must be (open) discs. Hence we obtain a CW-structure on the surface
  with $1$-cells $e_1, \ldots, e_k$ that are all loops. In particular,
  the $e_k$ generate the first homology of the surface.  Since the
  latter is nontrivial, at least one of the $e_k$ is a non-separating
  loop (as separating simple loops are null-homologous).
\end{proof}

Now we are ready to prove Lemma \ref{lem:fellowtravel}. The idea is that if the
width of the curve $\beta$ is large, then some vertical first
return must be fast. This in turn would yield a bounded size bicorn
beween the vertical foliation and $\beta$ -- which violates the Gromov
product bound.
\begin{proof}[{Proof of Small Width Lemma~\ref{lem:fellowtravel}}]
  We will show the counterpositive: that is, we show that if $\beta$
  contains a segment $\beta_0 \subset \beta$ so that the total
  horizontal width of $\beta_0$ is larger than $\epsilon$ but whose
  size is bounded by $B$, then the Gromov product
  $(\beta\cdot\xi_{\mathcal{F}_v})_{\beta_0}$ can be bounded from above.
  

  We let $L=L(B,\epsilon)$ be the bound from
  Lemma~\ref{lem:target-lemma}. Hence, any vertical geodesic of length
  at least $L$ for the singular flat metric determined by $q$
  intersects the segment $\beta_0$.
  
  By Lemma~\ref{lem:bicorn-existence}, there is a choice of such a
  (vertical) first return geodesic $\lambda_0$ (of length at most $L$)
  which defines a non-separating bicorn $\gamma$ with $\beta$. By
  construction, the size of this bicorn is bounded by $L+B$, and hence
  by Lemma~\ref{lem:size-bounds-distance} so is the distance
  $d(\gamma,b) \leq R=R(L+B)$. 

  \smallskip We now choose a sequence of curves
  $C_i = \lambda_i \cup \tau_i$ formed by an increasing sequence of
  vertical geodesics
  $\lambda_0 \subset \lambda_1 \subset \lambda_2 \subset \cdots$ of
  lengths diverging to $\infty$ and (short) horizontal geodesics
  $\tau_i$. By Proposition~\ref{prop:boundary-from-foliation}, these
  converge to the boundary point $\xi_{\mathcal{F}_v}$.  For $i>1$,
  the curve $\gamma$ will be a bicorn defined by the curves $C_i$ and
  $\beta$ (as all $C_i$ contain the segment $\lambda_0$).

  By Lemma~\ref{lem:bicorns}, $\gamma$ is contained in a
  $K$--quasigeodesic connecting $\beta$ to $C_i$ (for any $i$). Since
  the distance of $\gamma$ to the base point $b$ was bounded above by $R=R(L+B)$,
  this implies that the Gromov product $(\beta \cdot C_i)_{\beta_0}$ is
  bounded from above for all $i$ by a bound depending only on $K$
  (which is independent of everything) and $L, R$ (which depend only
  on $B, \epsilon$).
\end{proof}

With this lemma in hand, we can now prove point (3) of our main result Theorem \ref{thm:foliation-convergence}. Recall that we call a connected subset of a manifold a {\em pseudo-leaf} if it is a {\em finite} union of singular leaves of a singular foliation $\mathcal{F}$ (with isolated singularities). 
\begin{proposition}
  Suppose that a sequence of curves $\alpha_i$ converges to the
  boundary point $\xi_{\mathcal{F}}$ associated to an ending foliation
  in the Gromov boundary of $\ncd(S)$. Then the size of the $\alpha_i$
  converges to infinity, and convergent subsequences of segments of
  $\alpha_i$ limit to subsets of either a leaf or pseudo-leaf in the
  Hausdorff sense. More precisely, for any $\epsilon,B$ and
  subintervals $J_i \subset \alpha_i$ of size at most $B$, any
  accumulation point of the $J_i$ in the Hausdorff topology is a
  (pseudo-)leaf segment of $\mathcal{F}$.
\end{proposition}

\begin{proof}
  By Lemma~\ref{lem:size-bounds-distance}, the size of the $\alpha_i$
  needs to diverge, since otherwise the distance to the basepoint
  would be bounded.
  
  Since the $\alpha_i$ converge to $\xi_{\mathcal{F}}$, applying
  the Small Width Lemma~\ref{lem:fellowtravel} we see that for any $\epsilon, L$ and
  large enough $i$, any length--$L$ subsegment is contained in an
  $\epsilon$--neighbourhood of a (pseudo-)leaf.

  This implies that any accumulation point of any segments (of bounded
  length) of the $\alpha_i$ are contained in (pseudo-)leaf segments.
\end{proof}

Finally we conclude the proof of Theorem \ref{thm:foliation-convergence} by characterizing stabilizers in the following:
\begin{proposition}\label{prop:fixed-points}
  Suppose that $F:S \to S$ is a homeomorphism which preserves the
  boundary point $\xi_{\mathcal{F}} \in \partial_\infty \cd(S) $ defined by a vertical foliation $\mathcal{F}_v$ of a BSF-Structure that is ending. Then $F$ preserves the foliation $\mathcal{F}_v$.
\end{proposition}
\begin{proof}
  First observe that there is a number $K>0$ so that the following
  holds: if $\sigma$ is any rectifiable segment of length at most $L$, then
  a lift of $F\sigma$ to the universal cover (of the $P$-resolving cover) has
  diameter at most $KL$. Namely, we can cover $S$ with a finite number
  $B_1, \ldots, B_k$ of balls, so that the lift of any path $\sigma$
  of length $\leq L$ can be covered with $cL$ lifts of the $B_i$ (for
  some constant $c$). Now, a lift of $F$ expands each lift of a ball
  to a set of bounded diameter, hence the claim.

  \smallskip To show the proposition, it suffices to show that any regular leaf segment $\sigma$ of $\mathcal{F}$ is mapped into a
  leaf of $\mathcal{F}$. To this end, choose a sequence of curves
  $C_i$ obtained by closing very long leaf segments of $\mathcal{F}$
  containing $\sigma$ by short horizontal segments. The curves $C_i$
  converge to $\xi_{\mathcal{F}_v}$ by Corollary~\ref{cor:convergence-crit}.
  Note that the image $F(\sigma)$ has size
  at most $K\, l(\sigma)$ and at least $l(\sigma)/K$ by the claim above (applied to $F$ and $F^{-1}$).

  On the other hand, the curves $F(C_i)$ also converge (as vertices in
  $\ncd(S)$) to $\xi_\mathcal{F}$, since $F$ preserves this boundary point
  by assumption. Hence, the Gromov product with $\xi_\mathcal{F}$
  diverges. Applying the Small Width Lemma~\ref{lem:fellowtravel} to any $\epsilon$
  and $B=Kl(\sigma)$ we see that $F(\sigma)$ is eventually contained
  in an $\epsilon$--neighbourhood of a (pseudo-)leaf. Since this is true for
  any $\epsilon$, we obtain that $F(\sigma)$ is contained in a (pseudo-)leaf. This leaf must be regular, since singular points cannot be locally equivalent to regular points under leaf preserving maps and the claim follows.
\end{proof}
Using completely analogous arguments, we can also show
\begin{proposition}
  Suppose that $\lambda$ is an ending filling geodesic lamination for
  some hyperbolic metric on $S\setminus P$, for $P$ finite.

  Then there is a point $\xi_\lambda \in \partial_\infty\cd(S)$
  defined by $\lambda$. If a sequence of curves $\alpha_i$ converges
  to $\xi_\lambda$, then the $\alpha_i$ converge to leaf segments
  uniformly on bounded size sets.
  
  If a homeomorphism $F$ fixes $\xi_\lambda$ (as a boundary point),
  then $F$ preserves $\lambda$ (as a lamination).
\end{proposition}
\begin{proof} We only briefly indicate which arguments
  need to be adapted.  The definitions of size and width, as well as
  the Target Lemma for laminations were discussed at the end of
  Section~\ref{sec:target}.

  To define the boundary point, we now argue exactly as in
  Proposition~\ref{prop:boundary-from-foliation}, with $L$ being leaf
  segments of the lamination. The proof of
  Corollary~\ref{cor:convergence-crit} then also works almost exactly
  the same: instead of closing long vertical segments of the BSF
  structure by short horizontal segments, we close long leaf segments
  of $\lambda$ with short geodesic segments.

  A version of Lemma~\ref{lem:bicorn-existence}, replacing the
  vertical flow by the geodesic flow along the lamination $\lambda$ is
  true by the same argument. Finally the proofs of the Small Width
  Lemma~\ref{lem:fellowtravel} and Proposition~\ref{prop:fixed-points}
  go through {\em mutatis mutandis}.
\end{proof}

We emphasise that the stabilisers of $\xi_\lambda$ for a lamination
are very large -- any homeomorphism supported on a complementary
region of $\lambda$ will fix $\xi_\lambda$.

\section{Applications to pseudo-Anosovs}
In this section, we collect some applications of the previous results
to the action of point-pushing pseudo-Anosovs on the fine curve graph.

We now consider a point-pushing pseudo-Anosov mapping class relative to a
finite set of points $P$, which is not quasi-invertible  as a
mapping class  of $S \setminus P$ in the sense of Bestvina-Fujiwara. Note that many such pseudo-Anosovs
exist (cf.\ Theorem \ref{thm:BF1}). We denote by $F$ a dynamical Thurston representative of  such a pseudo-Anosov mapping class and by $\mathcal{F}_h, \mathcal{F}_v$ the horizontal and
vertical foliations (which then determine the fixed points at infinity
for the action of $F$ on the fine curve graph).

First, we observe the following corollary of
Proposition~\ref{prop:fixed-points}.
\begin{corollary}
  Any homeomorphism $G$ which commutes with $F$ (or conjugates $F$ to a
  positive power of itself) preserves the foliations
  $\mathcal{F}_h, \mathcal{F}_v$. In particular,
  $G$ preserves the set $P_1 \subseteq P$ of angle--$\pi$--singularities (1-prongs).

  If the mapping class on $S\setminus P$ induced by $G$ is trivial,
  then $G$ is the identity. In particular, the natural map from the
  centralizer of $F$ (in the homeomorphism group) to the centralizer
  of $[F] \in Mcg(S \setminus P)$ is an isomorphism; and so the
  centralizer of $F$ is virtually cyclic.
\end{corollary}
\begin{proof}
  The first claim is immediate from
  Proposition~\ref{prop:fixed-points} (or simply the fact that
  $\mathcal{F}_v, \mathcal{F}_h$ are determined by the dynamics of
  $F$).  For the second, we consider the universal cover
  $\mathbb{H}^2$ of the $P_1$-resolving cover of $S$, and a lift
  $\widetilde{G}$ of $G$ to $\mathbb{H}^2$.  Since $G$ preserves
  $\mathcal{F}_h, \mathcal{F}_v$, the lift $\widetilde{G}$ preserves
  the lifts of these foliations.
  
  Then, if $G$ defines the trivial mapping class of $S\setminus P_1$,
  and therefore lifts to a trivial mapping class of the $P_1$-resolving cover,
  the lift $\widetilde{G}$ acts trivially on the boundary at infinity
  of $\mathbb{H}^2$. Note that no two regular leaves of
  $\widetilde{\mathcal{F}}_h$ can share both endpoints, as otherwise there
  would be a foliated strip with isolated leaves contradicting
  minimality. The same is true for $\widetilde{\mathcal{F}}_v$. This
  implies that $\widetilde{G}$ actually fixes each individual leaf of
  $\widetilde{\mathcal{F}}_h, \widetilde{\mathcal{F}}_v$ setwise.  A
  dense set of points $p \in \mathbb{H}^2$ has the property that there
  is exactly one regular leaf of each of
  $\widetilde{\mathcal{F}}_h, \widetilde{\mathcal{F}}_v$ through $p$.
  Hence, $\widetilde{G}$ needs to fix all those $p$, and is therefore
  the identity homeomorphism and then so is $G$. 
\end{proof}
We next come to one of our main results, which gives the first explicit examples of elements in $\mathrm{Homeo}(S)$ that have positive stable commutator length. Moreover, this even holds if one considers commutators in the full group of orientation preserving homeomorphisms, generalizing our results in \cite{Dagger1}.
\begin{theorem}\label{thm:scl}
  Let $F$ be the Thurston representative  of point-pushing pseudo-Anosov mapping class relative to its
set of $1$-prongs  $P=P_1$ so that $[F] \in \mathrm{Mcg}(S\setminus P)$ is not
  conjugate to its inverse. Then $F$ has positive $\textrm{scl}(F)$ in
  $\mathrm{Homeo}_+(S)$.
\end{theorem}
\noindent To prove this theorem, we will show that a quasi-axis for $F$ acting on $\ncd(S)$ cannot be
quasi-inverted. By the main result of Bestvina-Fujiwara \cite{BF} this will imply the
proposition (cf.\ Section \ref{sec:BF}). The main technical work lies in the following lemma.

\begin{lemma}\label{lem:controlling-pi}
  Let $q=(\mathcal{F}_h, \mathcal{F}_v, d)$ be a BSF structure, and
  $P$ the set of angle-$\pi$ singularities. Assume that the horizontal
  and vertical foliations $\mathcal{F}_h, \mathcal{F}_v$ are
  ending. For any $\epsilon>0$ there is a $K=K(\varepsilon)>0$ so that
  the following holds. If $f$ is a homeomorphism so that
  \[ (f(\xi_{\mathcal{F}_v})\cdot\xi_{\mathcal{F}_h})_{\beta_0} > K, \quad
    (f(\xi_{\mathcal{F}_h})\cdot\xi_{\mathcal{F}_v})_{\beta_0} > K \] then
  $d(f(p), p) < \epsilon$ for any $p\in P$. 
\end{lemma}

\begin{figure}[h!]\label{fig:bell}
  \centering
  \includegraphics[height=0.5\textwidth]{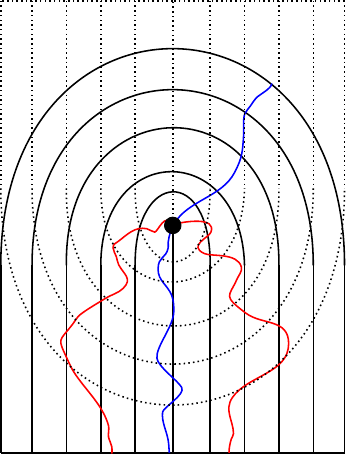}
  \caption{The ``bells'' in the proof of
    Lemma~\ref{lem:controlling-pi}. The solid black foliation is the
    vertical one, bending around the angle--$\pi$ singularity; the
    horizontal foliation is drawn dashed. In red, the image of the
    horizontal foliation under the homeomorphism is drawn. It cannot
    cross the vertical foliation too much (without bounding the Gromov
    product), and so it needs to exit at the bottom. But then, the
    (blue) image of the vertical foliation will have to also exit at
    the bottom, crossing the horizontal (dashed) foliation. }
  \label{fig:bell}
\end{figure}
\begin{proof}
  We argue the contrapositive. Namely, we will construct a
  neighbourhood $V = V_P$ of the angle-$\pi$ singularities, and show
  that if the Hausdorff distance satisfies $d(f(P),V) \ge \epsilon $, then there
  is a $K=K(\epsilon)>0$ so that one of the claimed inequalities
  fails.

  The subset $V$ consists of small neighbourhoods of each of the
  points $p_i \in P$ and each of these will in turn
  consist of the union of two ``bells'' $U_i, U_i'$. We begin by
  describing one of them. Denote by $\nu$ the vertical leaf segment of
  length $\epsilon$ starting at a point $p_i$ in the set $P$,
  and by $\eta$ the horizontal leaf extending length $\epsilon/2$ in
  both directions through the regular endpoint of $\nu$. We let $U_i$
  be the region swept out by the vertical arcs starting and ending on
  $\eta_i$ and bounding a disk containing $p_i$ together with a
  subinterval of $\eta_i$ (see Figure~\ref{fig:bell}). The boundary of
  $U_i$ consists of the \emph{base arc} $\eta_i$, and a \emph{top
    arc} lying on some leaf. We define the other bell $U'_i$ analogously, reversing the
  roles of vertical and horizontal.

  \smallskip The strategy of the proof is as follows: we will argue
  that leaves of the image of the vertical foliation
  $f(\mathcal{F}_v)$ need to leave the bell $U_i$ through the base,
  and similarly leaves of $f(\mathcal{F}_h)$ leave $U'_i$ through its
  base. This will lead to a contradiction to the assumption that no $1$-prong of $f(\mathcal{F}_v)$ lies in the union of the bells, since then there are two
  leaves of $f(\mathcal{F}_v)$ that cross transversely at a single point, forcing one to have large vertical or horizontal width.

  \smallskip We now give more details: Assume that for some point $p \in P$,
  we have $f(p) \notin V$, then (by a counting argument), $f(P)$
  is disjoint from some component $U_i \cup U'_i$ of $V$.

  Now observe that by construction, any segment $\sigma$ contained in
  $U_i$ has size at most $2\epsilon$. If such a segment $\sigma$ joins
  $p_i$ to the top arc, then $\sigma$ has total width $\epsilon$.
  Consider now the foliation $f(\mathcal{F}_h) \cap U_i$ on $U_i$, and
  observe that it does not have an angle-$\pi$-singularity. In
  particular, any leaf segment can be continued until it leaves the
  bell $U_i$ (since the only way leaf segments in a singular foliation
  cannot be continued is if they reach an angle-$\pi$ singularity, and
  leaves of minimal foliations eventually leave every disk).

  Let $\sigma$ be any maximal leaf segment of $f(\mathcal{F}_h) \cap
  U_i$ through $p_i$. If
  $\sigma$ had an endpoint on the top arc of
  $U_i$, it would be an arc of total width
  $\geq\epsilon$ and size at most $2\epsilon$. Let
  $\alpha_i$ be a sequence of curves formed by closing long horizontal
  leaf segments containing
  $f^{-1}(\sigma)$, which defines the boundary point
  $\mathcal{F}_h$ (as in Corollary~\ref{cor:convergence-crit}). Then
  for all large $i$, the curves $f(\alpha_i)$ contain
  $\sigma$ and, by the Small Width Lemma~\ref{lem:fellowtravel}, this yields an upper
  bound for the Gromov product $(f(\alpha_i)\cdot \mathcal{F}_v)_{\beta_0} <
  K(\epsilon)$. By definition, we then also get
  $(f(\mathcal{F}_h))\cdot \mathcal{F}_v)_{\beta_0} <
  K(\epsilon)$. By choosing
  $K > K(\varepsilon)$, this is impossible.  Hence, $\sigma$ leaves
  $U_i$ through the base.

  An analogous argument with the roles of horizontal and vertical reversed shows
  that any maximal leaf segment $\tau$ of $f(\mathcal{F}_v) \cap U'_i$
  through $p$ leaves $U'_i$ through its base.

  \smallskip But, since we assumed $f(P) \cap V_i= \emptyset$, we can
  choose leaf segments $\sigma, \tau$ as above, which only intersect
  in $p_i$ and cross there transversely -- since
  $f(\mathcal{F}_v) \cap (U'_i\cup U_i)$ and
  $f(\mathcal{F}_h) \cap (U'_i\cup U_i)$ are singular foliations on a
  disk without angle-$\pi$ singularities, and therefore such leaf segments
   cannot bound bigons.

  This contradicts the fact that both $\sigma, \tau$ leaf segments of $f(\mathcal{F}_h),f(\mathcal{F}_v)$ respectively must leave their bells
  through the bottom since any embedded arc disconnects a disk
  (compare Figure~\ref{fig:bell}).
\end{proof}

We are now ready to prove Theorem~\ref{thm:scl}.
\begin{proof}[Proof of Theorem~\ref{thm:scl}]
Let $C_0$ be any smooth non-separating simple closed curve disjoint from the punctures $P \subset S$ of the pseudo-Anosov map $F$ and set $C_i = F^i(C_0)$. Then $(C_i)$ is a quasi-geodesic axis of quality $L$ for the action of $F$ on $\ncd(S)$ and considering isotopy classes $([C_i])$ we similarly have a quasi-geodesic axis for the action of the mapping class on $\mathcal{NC}(S\setminus P)$ of quality $L' \ge L$.

Suppose that there is a sequence $g_k$ of homeomorphisms so
  that   as points in $\partial_{\infty} \ncd(S)$
  $$g_k\mathcal{F}_h \to \mathcal{F}_v, g_k\mathcal{F}_v \to
  \mathcal{F}_h.$$
By hyperbolicity of $\ncd(S)$, the following holds: for any $N$ and any $k$ large enough, 
  we have that if $i \geq N'(N)>N$, then $g_k(C_i)$ is within distance $B$ of
  $C_{-j}$ for some $j > N$ and, analogously, $g(C_{-i})$ is within
  distance $B$ of $C_j$ for some $j > N$. Here, $B$ only depends on
  the hyperbolicity constant and the quality $L$ of the quasigeodesic.

  \smallskip 
  The quasi-axis $([C_i])$ in $\mathcal{NC}(S\setminus P)$ is
  not quasi-invertible by our assumption on the mapping class $[F]$. In particular, there are $N, K$ so that there
  is no mapping class $\varphi$ of $S\setminus P$ so that
  \[ (\varphi C_i \cdot C_{-j})_{\beta_0} > K - B, \quad (\varphi C_{-i}\cdot C_j)_{\beta_0} > K - B. \]
  for all $i, j \geq N$ 

  \smallskip Now, let $\epsilon$ be such that the $2\epsilon$--neighbourhood of the puncture set $P$ is
  disjoint from $C_{i}$ for all $|i| \le N,N'$. Then by Lemma~\ref{lem:controlling-pi}, for $k$ large enough, we have that $g_k(P)$ is within
  $\epsilon$ of $P$, and also $g_k(C_{N'}), g_k(C_{-N'})$ are within distance $B$
  of $C_{-i}, C_j$ for $i,j > N$.

  Now, we can change $g_i$ by a $C^0$-small isotopy supported in a $2\varepsilon$-ball about $P$ to obtain a homeomorphism $g'$ with
  $g'(C_{N'}) = g_i(C_{N'}), g'(C_{-N'}) = g_i(C_{-N'})$ which {\em preserves} $P$. Hence, $g'$
  defines a mapping class $\varphi$ which violates the assumption that the quasi-axis $([C_i])$ is not quasi-invertible.
\end{proof}
\begin{remark}\label{rem:single_push}
  Suppose $F$ is a Thurston representative of a point pushing
  pseudo-Anosov $\psi$ of a single point $p \in S$. Then $F$ has positive
  stable commutator length in $\mathrm{Homeo}_0(S)$ without any
  further assumptions.

  Namely, in the above proof, the $g_i$ are then isotopic to the
  identity, and therefore so is $g'$. Hence, $g'$ is a representative
  of a point push (of the same point $p$).

  Applying \cite{BF} to the action of the point-pushing group
  $\pi_1(S,p)$ on the curve graph $\mathcal{C}(S-p)$ this would imply
  that $\psi$ is conjugate to its inverse in $\pi_1(S,p)$.  However,
  this is impossible since surface groups are bi-orderable.
\end{remark}

Finally, we prove a kind of Tits alternative which generalizes \cite[Theorem 1.10]{HurtadoXue} to all regularities and genera.
\begin{theorem}
  Suppose that $G < \mathrm{Homeo}_+(S)$ is a group which contains a Thurston representative
  of a pseudo-Anosov element of $\mathrm{Mcg}(S\setminus P)$ for some finite $P$.

  Then either $G$ contains a free group, or $G$ has an index $2$ subgroup which fixes a
  foliation on $S$.
\end{theorem}
\begin{proof}
  Let $F$ be the Thurston representative, and let $\xi_+, \xi_-$ be
  the boundary points of $\cd(S)$ defined by the stable and unstable
  foliations.  There are three cases. If $G$ preserves the set
  $\{\xi_+, \xi_-\}$, then there is an index $2$ subgroup which fixes
  both, and we can apply Proposition \ref{prop:fixed-points}. If
  $\xi_+$ or $\xi_-$ are fixed by $G$, we are also done.  Otherwise,
  there is a conjugate of $F$ in $G$ which has fixed points that are
  both distinct from $\xi_-, \xi_+$. But then, large enough powers of
  $F$ and this conjugate generate a free group by the ping pong lemma
  for actions on hyperbolic spaces.
\end{proof}

\begin{remark}
  We remark that since the (un)stable foliation of a pseudo-Anosov is
  uniquely ergodic, if the the group $G$ in the previous theorem does
  not contain a free group, then it preserves the transverse measure
  of the foliation (up to scale) as well.

  This is reminicent of the measurable Tits alternative by Margulis
  \cite{Margulis} in dimension $1$: a subgroup of
  $\mathrm{Homeo}(S^1)$ contains a free group or preserves a
  probability measure.
\end{remark}
\section{Stabilisers of Nonorientable Foliations}
In this section we put restriction on possible fixed points of homeomorphisms acting parabolically. This answers a question posed to us by Yair Minsky.
In particular, we have the following
\begin{theorem}\label{thm:parabolics}
  Suppose that $\mathcal{F}$ is a singular foliation on $S$ induced by a BSF
  structure without compact pseudo-leaves and singular set $Q$.
  Then the corresponding Gromov boundary point $\xi_\mathcal{F}$ is not the
  fixed point of any homeomorphism acting parabolically.

  If $\xi_\mathcal{F}$ is the fixed point of a homeomorphism acting
  hyperbolically, then the homeomorphism fixes $Q$ and defines a
  pseudo-Anosov mapping class on $S\setminus Q$ so that $\mathcal{F}$ is one of
  its invariant foliations.
\end{theorem}
\begin{proof} We first show the statement about parabolics, and argue
  by contradiction. Suppose that $\phi$ is a homeomorphism fixing
  $\xi_\mathcal{F}$, and suppose that $\phi$ acts as a parabolic isometry on
  $\cd(S)$. By Proposition~\ref{prop:fixed-points},
  $\phi$ then preserves the foliation $\mathcal{F}$, and in particular all
   singularities in $Q$.
   
The induced mapping class $[\phi] \in \mathrm{Mcg}(S\setminus Q)$ fixes the
isotopy class $[\mathcal{F}]$ of the foliation $\mathcal{F}$ defined on $S\setminus Q$. We remark that the foliation is minimal on $S\setminus Q$ due to the assumption that it has not compact pseudo-leaves  -- see discussion after the definition of ending in Section \ref{subsec:bsf-structures}. The mapping
class $[\phi]$ is therefore of finite order: no reducible infinite order mapping
class can fix a minimal foliation at all, and if $[\phi]$ were
pseudo-Anosov, then $\mathcal{F}$ would act on $\mathcal{C}^\dagger(S)$ as a
hyperbolic isometry, which is not the case by assumption.

Hence, up to replacing $\phi$ by a power we may assume that $[\phi]$
is the trivial mapping class on $S \setminus Q$. In particular,
lifting to the universal cover $\widetilde{X}$ of the orientation
cover (i.e. the cover which is branched with order $2$ at each
singularity in $Q$), a lift $\widetilde{\phi}$ of $\phi$ maps each leaf of
$\widetilde{\mathcal{F}}$ to itself. For all non-regular leaves,
additionally all prongs can be assumed to be fixed pointwise.

\smallskip Consider now any point $p \in \widetilde{X}$, and let $L_p$
be either the regular leaf containing $p$ (if it is regular), or the
singular leaf containing $p$ (i.e. the maximal ray or segment
containing $p$ bounded by singular points). We choose an orientation
on $\widetilde{\mathcal{F}}$ so that for the induced identification of $L_p$
with $\mathbb{R}$ (or $[0,\infty)$ or $[a,b]$), the point
$\widetilde{F}(p)$ is to the right of $p$. We aim to show that there
is some point $q \in L_p$ to the right of $\widetilde{F}(p)$ whose
image $\widetilde{F}(q)$ is to the left of $q$ -- and therefore that
there is an interval $I_p \subset L_p$ containing $p$ which is mapped
into itself by $\widetilde{\mathcal{F}}$.

To see this, note first that (by continuity) any $p'$ close enough to
$p$ also moves to the right in its $L_{p'}$. On the other hand (by
minimality of $\mathcal{F}$), we can find such $p'$ which map on $S$
into the image of $L_p$, but (by non-orientability of $\mathcal{F}$)
so that the deck transformation $\phi_{p,p'}$ which sends $L_p$ to
$L_{p'}$ reverses the orientation on the leaf.  This means that
$\phi_{p,p'}^{-1}(p') = q \in L_p$ has the desired
property. 

\smallskip
As a consequence, the images $\widetilde{\phi}^n(K)$ of a fundamental
domain $K$ are contained in a compact set $C$ independent of $n$. This implies that, given any curve $\beta$, the size of
$\phi^n(\beta)$ is bounded -- which, by
Lemma~\ref{lem:size-bounds-distance}, implies that $\phi$ is elliptic.

\smallskip The claim on hyperbolically follows in an analogous way -- if
$\phi$ acts hyperbolically, but $[\phi]$ is not a pseudo-Anosov, then
again $\phi$ would fix every leaf setwise, allowing for the same
contradiction.
\end{proof}
We remark that in the other direction, if one begins with an orientable foliation, then this is the fixed point of a parabolic element given by flowing along leaves as described in \cite{dagger2}.

\section{Lifting under finite covers and Quasi-morphisms}
In this section we consider connected surfaces $\Sigma$ with non-empty boundary of negative Euler characteristic (possibly non-orientable). Then as $\Sigma$ is hyperbolic by Hamstrom \cite{Hamstrom} the identity component $\Homeo_0(\Sigma,\partial \Sigma)$ of the groups of homeomorphisms fixing the boundary (pointwise) is contractible. 
Now consider any finite cover $\overline{\Sigma} \to \Sigma$. Any isotopy to the identity fixing the boundary lifts under a finite cover and thus there is well-defined lifting map. Composing this with the map induced by collapsing the boundary components of $\overline{\Sigma}$ to points we obtain 
$$ \Homeo_0(\Sigma,\partial \Sigma) \stackrel{\textrm{lift}}{\longrightarrow} \Homeo_0(\overline{\Sigma},\partial \overline{\Sigma}) \stackrel{\textrm{collapse}}{\longrightarrow} \Homeo_0(\widehat{\Sigma},P_{\partial}) \stackrel{\textrm{forget}}{\longrightarrow} \Homeo_0(\widehat{\Sigma}).$$
Using this map and the fact that Thurston representatives on $\Sigma$ map to Thurston representatives on $\widehat{\Sigma}$ one can deduce the following:
\begin{theorem}\label{thm:cover}
Let $\Sigma$ be any hyperbolic surface (possibly non-orientable) but with non-empty boundary. Then the group $\Homeo_0(\Sigma,\partial \Sigma)$ contains elements with $\mathrm{scl} >0$.
\end{theorem}
\begin{proof}
  First consider a Thurston representative $\varphi_T$ of a point pushing
  pseudo-Anosov of $\Sigma$ relative to some point $P \in int(\Sigma)$ so that the mapping class becomes trivial in $Mcg(\Sigma, \partial \Sigma)$ after
  forgetting the point $P$. We then choose a finite cyclic cover that is non-trivial on all boundary components.  In particular, the lift of the Thurston representative can have $1$-prongs only on the pre-image $\overline{P}$ of the marked point. After collapsing each boundary component we then obtain a pseudo-Anosov map on the closed surface $\widehat{\Sigma}$ having $1$-prongs precisely at the puncture set $\overline{P}$.
  
  Choosing
  $\overline{\Sigma}$ appropriately, we can assume that the closed
  surface $\widehat{\Sigma}$ is hyperbolic. Hence, the resulting map
  $\widehat{\varphi}_T$ is again a Thurston representative of some
  point-pushing map with respect to its $1$-prong set that we identify with $\overline{P} \subseteq int(\overline{\Sigma})$. Choosing $\varphi_T$ correctly (using
  Theorem~\ref{thm:BF1} and Theorem~\ref{thm:scl}) we can thus ensure
  that it has positive scl in $\Homeo(\widehat{\Sigma})$. It follows
  that $\varphi_T$ also has positive scl in $\Homeo_0(\Sigma)$.

  We remark that the Thurston representative $\varphi_T$ does not fix the boundary pointwise and acts with source-sink dynamics with the number of sinks corresponding to the order of the prong at a particular boundary component. In particular $\varphi_T$ does not lie in $\Homeo_0(\Sigma,\partial \Sigma)$ so we have only shown that this element has positive scl in the group $\Homeo_0(\Sigma)$ of homeomorphisms that are isotopic to the identity, but do not necessarily fix the boundary.  To
  address this, we argue as follows.

The assumption that the map is trivial in the mapping class
  group of $\Sigma$ implies that the homeomorphism $\varphi_T$ has {\em fixed
    points} on each boundary component that correspond to prongs of the stable/unstable
  foliations of the pseudo-Anosov map. We then attach a map of annuli
  $A$ that are given by an isotopy of the boundary map to the
  identity to obtain a `fattened' version of the
  Thurston representative
  $F^\mathrm{fat}_T$. By choosing the isotopies to have the correct number of twists, we can ensure that $F^\mathrm{fat}_T \in \Homeo_0(\Sigma,\partial \Sigma)$ is isotopic to the identity relative to the
  boundary. Conjugating $F^\mathrm{fat}_T$ by maps ``compressing the
  annulus'' $A$ to a an annulus $A_\epsilon$ of width $\epsilon$, we
  obtain a family of homeomorphisms $F^\epsilon_T$ which are all
  conjugate to each other in $\Homeo_0(\Sigma,\partial
  \Sigma)$. Outside the annulus $A_\epsilon$, the map $F^\epsilon_T$
  agrees with $F_T$.

  The maps $F^\epsilon_T$ do not converge as $\epsilon\to 0$ on
  $\Sigma$, but their images in
  $\Homeo_0(\widehat{{\Sigma}},P_{\partial})$ under the lifting and
  collapsing maps above {\em do converge} to the Thurston
  representative $\widehat{\varphi}_T$ on $\widehat{{\Sigma}}$ (since
  the annuli $A_\epsilon$ lift and collapse to disks of radius
  $\epsilon$). Since stable commutator length is
  $\mathcal{C}^0$--continuous (cf.\ \cite[Theorem 1.5]{Dagger1}), this
  implies that for small $\epsilon$ the image of $F^\epsilon_T$ has
  positive stable commutator length (and, since they are all
  conjugate, in fact all of them do).
\end{proof}
\noindent This result is sharp for orientable surfaces as we know that in the case of the annulus and the disc the groups do not admit quasi-morphisms \cite {BIP} and it was shown in \cite{Dagger1} that if $\Sigma = T^2$ is a torus then there are quasi-morphisms. The remaining cases are then the Klein bottle, the M\"obius band and the real projective plane.

\bibliographystyle{alpha}
\bibliography{boundary}

\end{document}